\documentclass{amsart}

\usepackage{amssymb, amsmath}
\usepackage{mathrsfs}
\usepackage{amscd}
\usepackage{verbatim}

\allowdisplaybreaks
\theoremstyle{plain}
\newtheorem{theorem}{Theorem}[section]
\theoremstyle{remark}
\newtheorem{remark}[theorem]{Remark}
\newtheorem{example}[theorem]{Example}
\theoremstyle{plain}

\newtheorem{lemma}[theorem]{Lemma}
\newtheorem{proposition}[theorem]{Proposition}
\newtheorem{definition}[theorem]{Definition}

\numberwithin{equation}{section}

\def\R{{\mathbb R}}
\def\C{{\mathbb C}}

\newcommand{\E}{{\mathbb E}}
\renewcommand{\P}{{\mathbb P}}
\newcommand{\F}{{\mathcal F}}

\newcommand{\g}{\gamma}

\renewcommand{\O}{\Omega}

\newcommand{\Div}{\operatorname{div}}

\newcommand{\calL}{{\mathcal B}}
\newcommand{\B}{{\mathcal B}}
\newcommand{\n}{\Vert}

\newcommand{\lb}{\langle}
\newcommand{\rb}{\rangle}

\newcommand{\tX}{{\widetilde X}}

\begin{document}

\title
[Stochastic equations with boundary noise]{Stochastic equations with
boundary noise}

\author{Roland Schnaubelt}
\address{Institute of Analysis, Department of Mathematics \\
Karlsruhe Institute of Technology (KIT)\\
D-76128  Karls\-ruhe\\Germany}
\email{schnaubelt@math.uni-karlsruhe.de}

\author{Mark Veraar}
\address{Delft Institute of Applied Mathematics\\
Delft University of Technology \\ P.O. Box 5031\\ 2600 GA Delft\\The
Netherlands}
\email{m.c.veraar@tudelft.nl, mark@profsonline.nl}

\thanks{The second named author was supported by the Alexander von Humboldt
foundation and by a 'VENI subsidie' (639.031.930) in the
'Vernieuwingsimpuls' programme of the Netherlands Organization for
Scientific Research (NWO)}

\begin{abstract}
We study the wellposedness and pathwise regularity of semilinear
non-autonomous parabolic evolution equations with boundary and interior noise
in an $L^p$ setting. We obtain existence and uniqueness of mild and weak
solutions. The boundary noise term is reformulated as a perturbation of a
stochastic evolution equation with values in extrapolation spaces.
\end{abstract}

\subjclass{Primary 60H15; Secondary 35R60, 47D06}


\keywords{Parabolic stochastic evolution equation, multiplicative
boundary noise,  non-autonomous equations, mild solution,
variational solution,  extra\-polation}

\date\today

\maketitle

\section{Introduction}
In this paper we investigate the wellposedness and pathwise regularity
of semilinear non-autonomous parabolic
evolution equations with boundary noise.
A model example which fits in the class of problems we study is given by
\begin{align}\label{eq:SE0}
\frac{\partial u}{\partial t}(t,s)
  &= \mathcal{A}(t,s,D) u(t,s)
 \qquad\text{ on } (0,T]\times S, \notag \\
\mathcal{B}(t,s,D) u(t,s) &= c(t,u(t,s)) \frac{\partial w}{\partial t}(t,s)
\qquad\text{ on } (0,T]\times \partial S, \\
u(0,s) &= u_0(s), \qquad\text{ on } S. \notag
\end{align}
Here $S\subset \R^d$ is a bounded domain with $C^2$ boundary,
 $\mathcal{A}(t,\cdot,D)=\Div (a(t,\cdot)\nabla)$ for uniformly positive
 definite, symmetric matrices $a(t,s)$
with the conormal boundary operator $\mathcal{B}(t,s,D)$,
$c(t,\xi)$ is Lipschitz in $\xi\in \C$,
$(w(t))_{t\ge0}$ is a Brownian motion for an filtration $\{\F_t\}_{t\ge0}$
and with values in $L^r(\partial S)$ for some $r\ge 2$, and
$u_0$ is an $\F_0$-measurable initial value. Actually, we also allow for
lower order terms, interior noise, nonlocal nonlinearities,
and more general stochastic terms, see Section~\ref{sec:ex1}.

As a first step one has to give a precise meaning to the formal
boundary condition in \eqref{eq:SE0}. We present two solution concepts
for \eqref{eq:SE0} in  Section~\ref{sec:ex1}, namely a mild and a weak one,
which are shown to be equivalent. Our analysis is then based on the mild
version of  \eqref{eq:SE0},  which fits into the general
framework of \cite{Vna}  where parabolic non-autonomous
evolution equations in Banach spaces were treated. The results in
\cite{Vna} rely on the stochastic integration theory in certain classes
of Banach spaces (see \cite{Brz2,NVW1,Pi75}).
In order to use \cite{Vna}, the inhomogeneous boundary term is reformulated as
an additive  perturbation of a stochastic
evolution equation corresponding to homogeneous boundary conditions.
This perturbation maps into a so--called extrapolation space
for the realization
$A(t)$  of $\mathcal{A}(t,\cdot,D)$ in $L^p(S)$ with the boundary condition
$\mathcal{B}(t,\cdot,D) u=0$ (where $p\in [2,r]$).
  Such an approach was
developed for deterministic problems by Amann in  e.g.\ \cite{Amnonhom} and
\cite{Am}. We partly use somewhat different techniques taken from
\cite{MSchn}, see also  the references therein. For this reformulation,
one further needs the solution map of a corresponding elliptic
boundary  value problem with boundary data in  $L^r(\partial S)$ which
is the range space of the Brownian motion. Here we heavily rely
on the theory presented in \cite{Amnonhom}, see also the references therein.
We observe that in \cite{Amnonhom} a large class of elliptic systems was
studied. Accordingly, we could in fact allow for systems in \eqref{eq:SE0},
but we decided to restrict ourselves to the scalar case in order to
simplify the presentation.

We establish in Theorem~\ref{thm:mainexdm} the existence and
 uniqueness of a mild solution $u$ to \eqref{eq:SE0}. Such a solution
 is  a process $u:[0,T]\times \Omega\to L^p(S)$ where $(\Omega,P)$
 is the probability space for the Brownian motion. We further show
 that for a.e.\ fixed $\omega\in \Omega$ the path $t\mapsto u(t,\omega)$
is (H\"older) continuous with values in suitable interpolation spaces
between $L^p(S)$ and the domain of $A(t)$, provided that $u_0$ belongs
to a corresponding interpolation space a.s.. As a consequence,
the paths of $u$ belong to $C([0,T],L^q(S))$ for all $q<dp/(d-1)$.
At this point, we make use of the additional regularity provided by
the  $L^p$ approach to stochastic evolution equations.

In \cite{Mas} an autonomous version of \eqref{eq:SE0} has been studied
in a Hilbert space situation (i.e., $r=p=2$) employing related techniques.
However, in this paper only regularity in the mean and no pathwise
regularity has been treated. In \cite[\S 13.3]{DPZ2}, Da Prato and
Zabczyk  have also investigated boundary noise of Neumann type. They
deal with a specific situation where $a(t)=I$, the domain is a cube
and the noise acts on one face which allows more detailed results.
See also \cite{AZ}, \cite{DFT} and \cite{Sow}  for further
contributions to problems with boundary noise. As explained in
Remark~\ref{rem:neumann} we cannot treat
Dirichlet type boundary conditions due to our methods. In one space
dimension Dirichlet boundary noise has been considered in \cite{AB2}
in weighted $L^p$--spaces by completely different techniques.

In the next section, we first recall the necessary material about
parabolic deterministic evolution equations and about stochastic
integration. Then we study an abstract stochastic
evolution equation related to \eqref{eq:SE0} in Section~\ref{sec:abstracteq}.
Finally, in the last section we treat a more general version of
\eqref{eq:SE0} and discuss various examples concerning the stochastic terms.

\section{Preliminaries\label{sec:prel}}
We write $a \lesssim_K b$ if there exists a constant
$c$ only depending on $K$ such that $a\le c b.$ The relation
$a\eqsim_K b$ expresses that $a \lesssim_K b$ and $b \lesssim_K a$.
If it is clear  what is meant, we just write $a \lesssim b$ for convenience.
Throughout, $X$ denotes a Banach space, $X^*$ its dual, and $\B(X,Y)$
the space of linear bounded operators from $X$ into another Banach space $Y$.
If the spaces are real, everything below should be understood
for the complexification of the objects under consideration.
The complex interpolation space for  an interpolation couple
 $(X_1,X_2)$ of order $\eta\in (0,1)$ is designated by
$[X_1,X_2]_{\eta}$.  We refer  to
\cite{Tr1} for the relevant definitions and basic properties.

\subsection{Parabolic evolution families}\label{sec:parabolic}
We briefly discuss the approach to non--auto\-no\-mous parabolic evolution
equations developed by Acquistapace and Terreni, \cite{AT2}. For $w\in \R$ and
$\phi\in [0,\pi]$, set $\Sigma(\phi,w) = \{\lambda\in \C: |\text{arg}(z-w)|\leq
\phi\}$. A family $(A(t),D(A(t)))_{t\in [0,T]}$ satisfies the hypothesis (AT)
if the following two conditions hold, where $T>0$ is given.

\begin{enumerate}
\item[(AT1)]  \label{AT1}
$A(t)$ are densely defined, closed linear operators on a
 Banach space $ X $ and there are
constants $w\in \R$, $K \ge 0$, and $ \phi \in (\frac{\pi}{2},\pi)$
such that $\Sigma(\phi,w) \subset \varrho(A(t))$ and
\[ \| R(\lambda, A(t)) \| \le \frac{K}{1+ |\lambda-w|}\]
holds  for all $\lambda \in \Sigma(\phi,w)$ and $t\in [0,T]$.
\item[(AT2)]  \label{AT2}
There are constants $L \ge 0$ and $\mu, \nu \in (0,1]$ such that
$\mu +\nu >1$ and
\[ \| A_w(t)R(\lambda,
A_w(t))(A_w(t)^{-1}- A_w(s)^{-1})\| \le L |t-s|^\mu (|\lambda|+1)^{-\nu}  \]
holds for all $\lambda \in \Sigma(\phi,0)$ and $s,t \in [0,T]$, where
$A_w(t)=A(t)-w$.
\end{enumerate}

Condition (A1) just means sectoriality with  angle $\phi>\pi/2$ and uniform
constants, whereas (A2) says that
the resolvents satisfy a H\"older condition in stronger norms.
In fact, Acquistapace and Terreni have studied a somewhat weaker version
of (AT2) and allowed for non dense domains. Later on, we work on reflexive
Banach spaces, where sectorial operators  are automatically densely defined
so that we have included the density assumption in (AT1) for  simplicity.
 The conditions (AT) and several variants of them have intensively been
studied in the literature, where also many examples can be found,
see e.g.\ \cite{Acq,AT2, Am, Schn, Ya90}.
If (AT1) holds and the domains $D(A(t))$ are constant in time,
then the H\"older continuity of $A(\cdot)$ in
$\calL(D(A(0)),X)$ with
exponent $\eta$ implies (AT2) with $\mu = \eta$ and $\nu=1$
(see \cite[Section~7]{AT2}).

Let $\eta\in(0,1)$,  $\theta\in [0,1]$, and $t\in [0,T]$. Assume that (AT1)
holds. The fractional power  $(-A_w(t))^{-\theta}\in \mathcal{B}(X)$ is
defined by
\[(-A_w(t))^{-\theta} = \frac{1}{2\pi i} \int_{\Gamma}
              (w-\lambda)^{-\theta} R(\lambda, A(t)) \, d\lambda,\]
where the contour
$\Gamma = \{\lambda: \arg(\lambda-w) = \pm \phi\}$ is
orientated counter clockwise. The operator
$(w-A(t))^{\theta}$ is defined as the inverse of $(w-A(t))^{-\theta}$.
We will also use the complex interpolation space
\[X_\eta^t = [X, D(A(t))]_{\eta}\,.\]
Moreover, the \emph{extrapolation space} $X^t_{-\theta}$ is the completion
of $X$ with respect to the norm $\|x\|_{X^t_{-\theta}} =
\|(-A_w(t))^{-\theta}x\|$. Let $A_{-1}(t):X\to X_{-1}^t$ be the unique
continuous extension of $A(t)$ which is sectorial of the same type. Then
$(w-A_{-1}(t))^\alpha:X_{-\theta}^t\to X_{-\theta-\alpha}^t$ is an
isomorphism, where $0\le \theta\le\alpha+\theta\le 1$. If $X$ is reflexive,
then one can identify the dual space $(X_{-1}^t)^*$ with $D(A(t)^*)$
endowed with its graph norm and the adjoint operator $A_{-1}(t)^*$
with $A(t)^*\in\B(D(A(t)^*),X^*)$. We mostly write $A(t)$ instead of
$A_{-1}(t)$. See e.g.\ \cite{Am, MSchn} for  more details.

Under condition (AT), we consider the non-autonomous Cauchy problem
\begin{equation}\label{nCP}
\begin{aligned}
u'(t) &= A(t) u(t), \qquad t\in [s,T],
\\ u(s) &= x,
\end{aligned}
\end{equation}
for  given $x\in X$ and $s\in[0,T)$.
A function $u$ is a {\em classical solution}\index{classical
solution}\index{solution!classical} of \eqref{nCP} if $u\in C([s,T];X)
\cap C^1((s,T]; X)$, $u(t)\in D(A(t))$ for all $t\in (s,T]$,
$u(s) = x$, and $\frac{d u}{dt}(t) = A(t) u(t)$ for all $t\in (s,T]$.
The solution operators of \eqref{nCP} give rise to the following
definition.
A family of bounded operators $(P(t,s))_{0\leq s\leq t\leq T}$ on
$X$ is called a {\em strongly continuous evolution
family}\index{evolution family}\label{p:evolu} if
\begin{enumerate}
\item $P(s,s) = I$  for all $s\in [0,T]$,
\item $P(t,s) = P(t,r) P(r,s)$ for all $0\leq s\leq r\leq t\leq T$,
\item the map $\{(\tau,\sigma)\in [0,T]^2: \sigma\leq \tau\}\ni (t,s)
\to P(t,s)$ is strongly continuous.
\end{enumerate}

The next theorem says that the  operators $A(t)$, $0\le t\le T$,
`generate' an evolution family having parabolic regularity. It is
a consequence of \cite[Theorem~2.3]{Acq}, see also \cite{AT2, Am, Schn, Ya90}.

\begin{theorem} \label{thm:exist-parab}
If condition (AT) holds, then there exists a unique strongly
continuous evolution family $(P(t,s))_{0\le s \le t \le T}$ such that
$u=P(\cdot, s) x$ is the unique classical solution of \eqref{nCP}
for every $x\in X$ and $s\in [0,T)$. Moreover, $(P(t,s))_{0\le s \le
t \le T}$ is continuous in $\calL(X)$ on $0 \le s < t \le T$ and
there exists a constant $C>0$ such that
\begin{align}\label{eq:2_13}
\|P(t,s) x\|_{X_{\alpha}^t} & \le  C(t-s)^{
\beta-\alpha}\|x\|_{X^{s}_{\beta}}
\end{align}
for all $0\leq\beta\leq\alpha \le 1$ and $0\leq s<t\leq T$.
\end{theorem}

\noindent
We further recall from \cite[Theorem 2.1]{Ya} that there is a constant $C>0$
such that
\begin{equation}
\label{eq:2_15-0} \| P(t,s)( w-A(s))^\theta x\|  \le  C (\mu-\theta)^{-1}(t-s)^{- \theta} \|x\|
\end{equation}
for all $0 \le s<t\le T$, $\theta\in (0,\mu)$ and $x\in
D((w-A(s))^{\theta})$. Clearly, \eqref{eq:2_15-0} allows to extend $P(t,s)$
to a bounded operator $P_{-\theta}(t,s):X_{-\theta}^s\to X$ satisfying
\begin{equation} \label{eq:2_15}
\| P_{-\theta}(t,s)( w-A_{-1}(s))^\theta\| \le
     C(\mu- \theta)^{-1}(t-s)^{- \theta}
\end{equation}
for all $0 \le s<t\le T$ and $\theta\in (0,\mu)$. Again, we mostly omit the
index $-\theta$.

\subsection{Stochastic integration\label{sec:stochint}}
Let $H$ be a separable Hilbert space with scalar product $[\cdot,\cdot]$,
$X$ be a Banach space, and
 $(S,\Sigma,\mu)$ be a measure space. A function $\phi:S\to X$ is called
 {\em strongly measurable} if it is the pointwise limit of a sequence of simple
functions. Let $X_1$ and $X_2$ be Banach spaces. An operator-valued function
$\Phi:S\to \calL(X_1,X_2)$ will be called {\em $X_1$-strongly measurable} if
the $X_2$-valued function $\Phi x$ is strongly measurable for all $x\in X_1$.

Throughout this paper  $({\O},{\F},{\P})$ is a probability space with a
filtration $(\F_t)_{t\geq 0}$ and $(\g_n)_{n\ge 1}$ is a {\em Gaussian
sequence}; i.e., a sequence of independent, standard,
real-valued Gaussian random variables defined on  $({\O},{\F},{\P})$.
An operator $R\in\calL(H,X)$ is said to be a {\em $\g$-radonifying
operator} if there exists an orthonormal basis $(h_n)_{n\ge 1}$ of $H$ such
that $\sum_{n\geq 1} \g_n \, Rh_n$ converges in $L^2(\O;X)$, see
\cite{Bog,DJT}. In this case we define
$$ \n R\n_{\g(H,X)}
:= \Bigl(\E\Bigl\n  \sum_{n\geq 1} \g_n \,Rh_n \Bigr\n^2\Bigr)^\frac12.$$ This
number does not depend on the sequence $(\g_n)_{n\ge 1}$ and the basis
$(h_n)_{n\ge 1}$, and defines a norm on the space $\gamma(H, X)$ of all
$\gamma$-radonifying operators from $H$ into $X$. Endowed with this norm,
$\g(H,X)$ is a Banach space, and it
holds $\n R\n\le \n R\n_{\g(H,X)}$. Moreover,
$\g(\mathcal{H},X)$ is an operator ideal in the
sense that if $S_1\in \calL(\tilde{\mathcal{H}}, \mathcal{H})$ and
$ S_2\in\calL(X,\tilde{X})$, then $R\in \g(\mathcal{H},X)$ implies
$S_2RS_1\in \g(\tilde{\mathcal{H}},\tilde{X})$ and
\begin{equation}\label{ideal}
 \n S_2RS_1\n_{ \g(\tilde{\mathcal{H}},\tilde{X})} \leq \n S_2\n \n R\n_{\g(\mathcal{H},X)} \n S_1\n.
\end{equation}

If $X$ is a Hilbert space, then $\g(H,X) = \mathcal{C}^2(H,X)$
isometrically, where $\mathcal{C}^2(H,X)$ is the space of
Hilbert-Schmidt operators. Also for $X=L^p$ there is a convenient
characterization of $R\in \g(H,L^p)$ given in \cite[Theorem 2.3]{BrzvN03}.
We use a slightly different formulation taken from \cite[Lemma~2.1]{NVW3}.

\begin{lemma}\label{lem:sq-fc-Lp}
Let $(S,\Sigma,\mu)$ be a $\sigma$-finite measure space and let $1\le
p<\infty$. For an operator $R\in\calL(H,L^p(S))$ the following assertions are
equivalent.
\begin{enumerate}
\item  $R\in \g(H,L^p(S))$.
\item  There exists a function $g\in L^p(S)$ such that for all $h\in H$ we have
$|R h|\le \|h\|_H\cdot g$ $\mu$-almost everywhere.
\end{enumerate}
Moreover, in this situation we have
\begin{equation}\label{eq:Rsquarekappa}
\|R\|_{\g(H,L^p(S))} \lesssim_p \|g\|_{L^p(S)}.
\end{equation}
\end{lemma}

A Banach space $X$ is said to have {\em type $2$} if there exists a constant
$C\ge 0$ such that for all finite subsets $\{x_1,\dots,x_N\}$ of $X$ we have
$$
\Big(\E \Big\n \sum_{n=1}^N r_n x_n\Big\n^2\Big)^\frac12 \le C
\Big(\sum_{n=1}^N \n x_n\n^2\Big)^\frac12.
$$
Hilbert spaces and  $L^p$-spaces with $p\in[2,\infty)$ have
type $2$. We refer to \cite{DJT} for details.
We will also need UMD Banach spaces. The definition of a UMD space will be
omitted, but we recall that every UMD space is reflexive. We refer to
\cite{Bu3} for an overview on the subject.  Important
examples of UMD spaces are the reflexive scale of $L^p$, Sobolev,
Bessel--potential and Besov spaces.

A detailed stochastic integration theory for operator-valued processes
$\Phi:[0,T]\times\O\to \calL(H,X)$, where $X$ is a UMD space, has been
developed in \cite{NVW1}. The full generality of this theory is not needed
here, since we can work with UMD spaces $X$ of type $2$ which allow for a
somewhat simpler theory.
Instead of of being a UMD space with type $2$, one can also assume
that $X$ is a space of
martingale type $2$ (cf.\ \cite{Brz2,Pi75}).

A family $W_H=(W_H(t))_{t\in \R_+}$ of bounded linear operators from $H$ to
$L^2(\O)$ is called an {\em $H$-cylindrical Brownian motion} if
\renewcommand{\labelenumi}{(\roman{enumi})}
\renewcommand{\theenumi}{(\roman{enumi})}
\begin{enumerate}
\item $\{W_H(t_j)h_k: j=1,\dots, J;\:  k=1, \dots,K\}$ is a Gaussian vector
 for all choices of $t_j\ge0$ and $h_k\in H$, and $\{W_H(t)h: t\geq
 0\}$ is a standard scalar Brownian motion with respect to the
 filtration $(\F_t)_{t\geq 0}$ for each $h\in H$;
\item $ \E (W_H(s)g \cdot W_H(t)h) = (s\wedge t)\,[g,h]_{H}$ for
all $s,t\in \R_+$ and  $g,h\in H.$
\end{enumerate}
\renewcommand{\labelenumi}{(\arabic{enumi})}
\renewcommand{\theenumi}{(\arabic{enumi})}

Now let $X$ be a UMD Banach space with type $2$. For an $H$-strongly measurable
and adapted $\Phi:[0,T]\times \O\to \g(H,X)$ which belongs to $L^2((0,T)\times
\O;\g(H,X))$ one can define the stochastic integral $\int_0^T \Phi(s) \, d
W_H(s)$ as a limit of integrals of adapted step processes, and
 there is a constant $C$ not depending on $\Phi$ such that
\[ \E\, \Big\|\int_0^T \Phi(s) \, d W_H(s)\Big\|^2
 \leq C^2 \|\Phi\|^2_{L^2((0,T)\times \O;\g(H,X))},\]
cf.\ \cite{Brz2}, \cite{NVW1}, and the references therein.
By a localization argument one may extend the class of integrable processes to
all $H$-strongly measurable and adapted $\Phi:[0,T]\times \O\to \g(H,X)$ which
are contained in $L^2(0,T;\g(H,X))$ a.s.. Below we use in particular
 the next result (see \cite{Brz2} and \cite[Corollary 3.10]{NVW1}).

\begin{proposition}\label{prop:stochinttype2}
Let $X$ be a UMD space with type $2$ and $W_H$ be a $H$-cylindrical
Brownian motion. Let $\Phi:[0,T]\times \O\to \g(H,X)$ be
$H$-strongly measurable and adapted. If $\Phi\in L^2(0,T;\g(H,X))$ a.s., then
$\Phi$ is stochastically integrable with respect to $W_H$ and for all $p\in (1,
\infty)$ it holds
\[ \Big(\E \sup_{t\in
[0,T]}\Big\|\int_0^{t}\Phi(s)\,dW_H(s)\Big\|^p\Big)^{\frac{1}{p}}
\lesssim_{X,p} \|\Phi\|_{L^p(\O;L^2(0,T;\g(H,X)))}.
\]
\end{proposition}

\noindent
In the setting of  Proposition~\ref{prop:stochinttype2} we also have,
for $x^*\in X^*$,
\begin{equation}\label{weak-int}
\Big\lb \int_0^T \Phi(s) \, d W_H(s),x^*\Big\rb
   = \int_0^T \Phi(s)^* x^* \, d W_H(s) \quad  \text{a.s.},
\end{equation}
cf.\ \cite[Theorem~5.9]{NVW1}.

\section{The abstract stochastic evolution equation\label{sec:abstracteq}}

Let $H_1$ and $H_2$ be separable Hilbert spaces, and let $X$ and $Y$ be
Banach spaces. On $X$ we consider the stochastic
evolution equation

\begin{equation}\tag{SE}\label{SEtype}
\left\{\begin{aligned}
dU(t) & = (A(t)U(t) + F(t,U(t)) + \Lambda_G(t) G(t,U(t)))\,dt  \\
& \qquad + B(t,U(t))\,dW_{H_1}(t) + \Lambda_C(t) C(t,U(t)) \, d W_{H_2}(t),
    \  t\in [0,T],\\
 U(0) & = u_0.
\end{aligned}
\right.
\end{equation}
Here $(A(t))_{t\in [0,T]}$ is a family of closed operators
on $X$ satisfying (AT). The processes $W_{H_1}$ and $W_{H_2}$ are independent
cylindrical Brownian motions with respect to $(\F_t)_{t\in [0,T]}$.
The initial value is a strongly $\F_0$-measurable mapping $u_0:\O\to X$.
We assume that the mappings $\Lambda_G(t): Y^t\to
X_{-\theta_G}^t$ and $\Lambda_C(t):Y^t\to X_{-\theta_C}^t$ are linear and
bounded, where the numbers $\theta_G,\theta_C\in[ 0,1]$ are specified below.
In Section \ref{sec:ex1}, the operators $\Lambda_G(t)$ and $\Lambda_C(t)$
are used to treat inhomogeneous boundary conditions. Concerning $A(t)$,
we make the following hypothesis.

\let\ALTERWERTA\theenumi
\let\ALTERWERTB\labelenumi
\def\theenumi{(H1)}
\def\labelenumi{(H1)}
\begin{enumerate}
\item \label{as:isom} Assume that $(A(t))_{t\in [0,T]}$ and
$(A(t)^*)_{t\in [0,T]}$ satisfy (AT) and that there exists an $\eta_0\in (0,1]$
and a family of Banach spaces $(\tX_{\eta})_{\eta\in [0,\eta_0]}$ such that
\[\tX_{\eta_0} \hookrightarrow  \tX_{\eta_1} \hookrightarrow  \tX_{\eta_2}
\hookrightarrow \tX_{0}=X \quad \text{ for all \ } \eta_0>\eta_1>\eta_2>0, \]
and each $\tX_{\eta}$ is a UMD space with type $2$. Moreover, it holds
\[[X, D(A(t))]_{\eta} \hookrightarrow \tX_{\eta}
  \quad \text{for all \ } \eta\in [0,\eta_0],\]
where the embeddings are bounded uniformly in $t\in[0,T]$.
\end{enumerate}
\let\theenumi\ALTERWERTA
\let\labelenumi\ALTERWERTB

\noindent
Assumption \ref{as:isom} has been employed in \cite{Vna} to deduce space time
regularity results for equations of the form \eqref{SEtype}, where spaces
such as $\tX_\eta$  have been used to get rid of the time dependence of
interpolation spaces; see also \cite[(H2)]{MSchnnew}. We have included an
assumption on $(A(t)^*)_{t\in [0,T]}$ for the treatment of
variational solutions. This could be done in a more general way as well,
but for us the above setting suffices. Assumption \ref{as:isom} can be
verified in many applications, see e.g.\ Section~\ref{sec:ex1}.

Let $a\in [0,\eta_0)$. The nonlinear terms $F,G,B$ and $C$ in \eqref{SEtype}
map as follows:
\[F:[0,T]\times\O\times \tX_a\to X, \qquad G(t):\O\times \tX_a\to Y^t,\]
\[B(t):\O\times \tX_a\to \g(H,X_{-1}^t), \qquad
  C(t):\O\times \tX_a\to \g(H_2,Y^t),\]
for each $t\in [0,T]$, where $Y^t$ are Banach spaces. We put
 $G(t)(\omega, x) = G(t,\omega, x)$, $B(t)(\omega,x) = B(t,\omega,x)$
and $C(t)(\omega,x) = C(t,\omega,x)$ for $(t,\omega,x)\in [0,T]\times \Omega \times X$.
Assuming \ref{as:isom} and $a\in [0,\eta_0)$, we state our main
hypotheses on $F,G,B$ and  $C$.

\let\ALTERWERTA\theenumi
\let\ALTERWERTB\labelenumi
\def\theenumi{(H2)}
\def\labelenumi{(H2)}
\begin{enumerate}
\item \label{as:LipschitzFtype} For all $x\in \tX_a$, the map
$(t, \omega)\mapsto F(t,\omega,x)$ is strongly measurable and adapted.
The function $F$ has linear growth and is Lipschitz continuous in space
uniformly in $[0,T]\times\O$; that is, there are constants $L_F$ and $C_F$
such that for all $t\in [0,T],$ $\omega\in \O$ and  $x,y\in \tX_{a}$ we have
\begin{align*}\nonumber
\|F(t,\omega,x)-F(t,\omega,y)\|_{X} & \leq  L_F\|x-y\|_{\tX_{a}},
\\ \nonumber
\|F(t,\omega,x)\|_{X}&\leq  C_F(1+\|x\|_{\tX_{a}}).
\end{align*}
\end{enumerate}
\let\theenumi\ALTERWERTA
\let\labelenumi\ALTERWERTB

\let\ALTERWERTA\theenumi
\let\ALTERWERTB\labelenumi
\def\theenumi{(H3)}
\def\labelenumi{(H3)}
\begin{enumerate}
\item \label{as:LipschitzGtype} For all $x\in \tX_a$, the map
$(t, \omega)\mapsto (-A_w(t))^{-\theta_G} \Lambda_G(t) G(t, \omega,x)\in X$
is strongly measurable
and adapted. The function $(-A_w)^{-\theta_G} \Lambda_G G$ has linear growth
and is Lipschitz continuous in space uniformly in $[0,T]\times\O$; i.e.,
there are constants $L_G$ and $C_G$ such that for all
$t\in [0,T],$ $\omega\in \O$ and  $x,y\in \tX_{a}$ we have
\begin{align*}\nonumber
\|(-A_w(t))^{-\theta_G} \Lambda_G(t)(G(t,\omega,x)-G(t,\omega,y))\|_{X}&\leq
L_G\|x-y\|_{\tX_{a}},
\\ \nonumber
\|(-A_w(t))^{-\theta_G} \Lambda_G(t) G(t,\omega,x)\|_{X}&\leq
C_G(1+\|x\|_{\tX_{a}}).
\end{align*}
\end{enumerate}
\let\theenumi\ALTERWERTA
\let\labelenumi\ALTERWERTB

\let\ALTERWERTA\theenumi
\let\ALTERWERTB\labelenumi
\def\theenumi{(H4)}
\def\labelenumi{(H4)}
\begin{enumerate}
\item \label{as:LipschitzBtype} Let $\theta_B\in [0,\mu)$ satisfy
$a+\theta_B <\frac12$. For all $x\in \tX_a$, the map $(t, \omega)\mapsto
(-A_w(t))^{-\theta} B(t, \omega,x)\in \g(H_1,X)$ is strongly measurable and
adapted. The function $(-A_w)^{-\theta_B}B$ has linear growth and is Lipschitz
continuous in space uniformly in $[0,T]\times\O$; that is, there are constants
$L_B$ and $C_B$ such that for all
$t\in [0,T],$ $\omega\in \O$ and  $x,y\in \tX_{a}$ we have
\begin{align*}
\nonumber
\|(-A_w(t))^{-\theta_B}(B(t,\omega,x)-B(t,\omega,y))\|_{\g(H_1,X)}&\leq
 L_B\|x-y\|_{\tX_{a}},
\\ \nonumber
\|(-A_w(t))^{-\theta_B}B(t,\omega,x)\|_{\g(H_1,X)}&\leq
C_B(1+\|x\|_{\tX_{a}}).
\end{align*}
\end{enumerate}
\let\theenumi\ALTERWERTA
\let\labelenumi\ALTERWERTB

\let\ALTERWERTA\theenumi
\let\ALTERWERTB\labelenumi
\def\theenumi{(H5)}
\def\labelenumi{(H5)}
\begin{enumerate}
\item \label{as:LipschitzCtype} Let $\theta_C\in [0,\mu)$ satisfy
$a+\theta_C<\frac12$. For all $x\in \tX_a$, the mapping  $(t, \omega)\mapsto
(-A_w(t))^{-\theta_C} \Lambda_C(t) C(t, \omega,x)\in\g(H_2,X)$ is strongly
measurable
and adapted. The function $(-A_w)^{-\theta_C} \Lambda_C C$ has linear growth
and is Lipschitz continuous in space uniformly in $[0,T]\times\O$; that is,
there are constants $L_G$ and $C_G$ such that for all
$t\in [0,T],$ $\omega\in \O$ and  $x,y\in \tX_{a}$ we have
\begin{align*}\nonumber
\|(-A_w(t))^{-\theta_C}
\Lambda_C(t)(C(t,\omega,x)-C(t,\omega,y))\|_{\g(H_2,X)}&\leq
L_C\|x-y\|_{\tX_{a}},
\\ \nonumber
\|(-A_w(t))^{-\theta_C} \Lambda_C(t) C(t,\omega,x)\|_{\g(H_2,X)}&\leq
C_C(1+\|x\|_{\tX_{a}}).
\end{align*}
\end{enumerate}
\let\theenumi\ALTERWERTA
\let\labelenumi\ALTERWERTB

We introduce our first solution concept.

\begin{definition}\label{def:mild}
Assume that \ref{as:isom}--\ref{as:LipschitzCtype} hold for
some $\theta_G,\theta_B,\theta_C\ge0$ and  $a\in  [0,\eta_0)$. Let
$r\in (2,\infty)$ satisfy $\min\{1-\theta_G, \frac12-\theta_B,
\frac12-\theta_C\}>\frac1r$. We call an $\tX_a$-valued process $(U(t))_{t\in
[0,T]}$ a {\em mild solution}\index{solution!mild}\index{mild solution} of
\eqref{SEtype} if
\renewcommand{\labelenumi}{(\roman{enumi})}
\renewcommand{\theenumi}{(\roman{enumi})}
\begin{enumerate}
\item $U:[0,T]\times\O\to \tX_a$ is strongly measurable and adapted,
and we have $U\in L^r(0,T;\tX_a)$ almost surely,

\item for all $t\in [0,T]$, we have
\[U(t) \!=\! P(t,0) u_0 + P*F(\cdot,U)(t)+P*\Lambda_G G(\cdot,U)(t) + P\diamond_1 B(\cdot,U)(t) +  P\diamond_2 \Lambda_C C(\cdot,U)(t)\]
in $X$ almost surely.
\end{enumerate}
\renewcommand{\labelenumi}{(\arabic{enumi})}
\renewcommand{\theenumi}{(\arabic{enumi})}
\end{definition}

Here we have used the abbreviations
\[P*\phi(t) = \int_0^t P(t,s) \phi(s) \, ds,\qquad
P\diamond_k \Phi(t) = \int_0^t P(t,s) \Phi(s) \, d W_{H_k}(s), \ \ k=1,2,\]
whenever the integrals are well-defined. Under our hypotheses
both $P*F(\cdot,U)(t)$ and
$P*\Lambda_G G(\cdot,U)(t)$ are in fact well-defined in $X$. Indeed, for the
first one this is clear from \ref{as:LipschitzFtype}. For the second
one we may write
\[P(t,s) \Lambda_G(s) G(s,U(s)) = P(t,s)(-A_w(s))^{\theta_G} (-A_w(s))^{-\theta_G} \Lambda_G(s) G(s,U(s))\]
It then follows from \eqref{eq:2_15}, H\"older's inequality,
and \ref{as:LipschitzGtype} that
\begin{align*}
\int_0^t \|P(t,s) & \Lambda_G(s) G(s,U(s))\|_{X} \, ds \\ & \lesssim \int_0^t
(t-s)^{-\theta_G} \|(-A_w(s))^{-\theta_G} \Lambda_G(s) G(s,U(s))\|_{X}  \, ds\\
& \lesssim 1+\|U\|_{L^{r}(0,T;\tX_a)},
\end{align*}
using that $1-\theta_G>\frac1r$. Similarly one can show that
$P\diamond_1 B(\cdot,U)(t)$ and $P\diamond_2 \Lambda_C
C(\cdot,U)(t)$ are well-defined in $X$, taking into account
Proposition~\ref{prop:stochinttype2}: Estimate \eqref{eq:2_15},
H\"older's inequality and \ref{as:LipschitzBtype} imply that
\begin{align*}
\int_0^t \|P(t,s)&  B(s,U(s))\|^2_{\g(H_1,X)} \, ds  \\&\lesssim \int_0^t
(t-s)^{-2\theta_B} \|(-A_w(s))^{-\theta_B} B(s,U(s))\|^2_{\g(H_1,X)}  \, ds\\ &
\lesssim 1+\|U\|_{L^r(0,T;\tX_a)}^2
\end{align*}
since $\frac12-\theta_B>\frac1r$. In the same way it can be proved  that
the  integral with respect to $W_{H_2}$ is well-defined.

We also recall the definition of a variational solution from \cite{Vna}.
To that purpose, for $t\in [0,T]$, we set
\begin{equation}\begin{aligned}\label{eq:Gammat}
\Gamma_t = \big\{ \varphi \in C^1([0,t];X^*) \,:
& \ \varphi(s)\in D(A(s)^*) \text{ for all } s \in [0,t]\\ & \text{ and }
[s\mapsto A(s)^*\varphi(s)] \in C([0,t];X^*) \big\}.
\end{aligned}
\end{equation}

\begin{definition}\label{def:var}
Assume that \ref{as:isom}--\ref{as:LipschitzCtype} hold with
$a\in  [0,\eta_0)$. An $\tX_a$-valued process $(U(t))_{t\in [0,T]}$ is
called a {\em variational solution}\index{variational
solution}\index{solution!variational} of \eqref{SEtype} if
\renewcommand{\labelenumi}{(\roman{enumi})}
\renewcommand{\theenumi}{(\roman{enumi})}
\begin{enumerate}
\item\label{en:var1} $U$ belongs to $L^2(0,T;\tX_a)$ a.s.\ and $U$ is strongly
measurable and adapted,

\item\label{en:var4} for all $t\in [0,T]$ and all $\varphi\in \Gamma_t$,
almost surely we have
\begin{align}\label{varsol}
\lb U(t), \varphi(t)\rb - \lb u_0, \varphi(0)\rb &=  \int_0^t [\lb U(s),
\varphi'(s) \rb + \lb U(s), A(s)^* \varphi(s) \rb \\
& \qquad + \lb F(s,U(s)), \varphi(s) \rb +
 \lb \Lambda_G(s) G(s,U(s)), \varphi(s) \rb]\,ds \notag\\
& \quad + \int_0^t B(s,U(s))^{*} \varphi(s) \, d W_{H_1}(s) \notag \\
& \quad + \int_0^t (\Lambda_C(s) C(s,U(s)))^{*} \varphi(s)\,dW_{H_2}(s).\notag
\end{align}

\end{enumerate}
\renewcommand{\labelenumi}{(\arabic{enumi})}
\renewcommand{\theenumi}{(\arabic{enumi})}
\end{definition}

\noindent
The integrand $B(s,U(s))^{*} \varphi(s)$  in \eqref{varsol} should be read as
\[((-A_w(s))^{-\theta_B} B(s,U(s)))^{*} (-A_w(s)^*)^{\theta_B}\varphi(s).\]
It follows from \ref{as:LipschitzBtype} that the function
$s\mapsto ((-A_w(s))^{-\theta_B} B(s,U(s)))^{*}$ is $X^*$-strongly measurable.
Moreover, the map
\[s\mapsto (-A_w(s)^*)^{\theta_B}\varphi(s)
  = (-A_w(s)^*)^{-1+\theta_B} (-A_w(s)^*) \varphi(s)\]
belongs to $C([0,t];X^*)$ by the H\"older continuity of
$s\mapsto (-A_w(s))^{-1+\theta_B}$ (cf. \cite[(2.10) and (2.11)]{Schn})
and the assumption on $\varphi$. Using \ref{as:LipschitzBtype}, we thus obtain
that the integrand is contained in $L^2(0,T;H_1)$
a.s.. As a result, the first stochastic integral in \eqref{varsol} is
well-defined. The other integrands have to be  interpreted similarly.

The next result shows that both solution concepts are equivalent in our
setting. It follows from Proposition~5.4 and Remark~5.3 in \cite{Vna} in the
same way as Theorem~\ref{thm:mainexistencegeninitial} below. (Remark~5.3 can
be used since $X$ is reflexive as a UMD space.)

\begin{proposition}\label{prop:varmild}
Assume that \ref{as:isom}--\ref{as:LipschitzCtype} hold for
some $\theta_G,\theta_B,\theta_C\ge0$ and  $a\in  [0,\eta_0)$. Let
$r\in (2,\infty)$ satisfy
 $\max\{\theta_C,\theta_B\}<\frac12-\frac1r$ and
$\theta_G<1-\frac1r$. Let $U:[0,T]\times\O\to \tX_a$ be a  strongly measurable
and adapted process such that  $U$ belongs to  $L^r(0,T;\tX_a)$ a.s..
Then $U$ is a mild
solution of \eqref{SEtype} if and only if $U$ is a variational solution of \eqref{SEtype}.
\end{proposition}

We can now state the main existence and regularity result for \eqref{SEtype}.

\begin{theorem}\label{thm:mainexistencegeninitial}
Assume that \ref{as:isom}--\ref{as:LipschitzCtype} hold for
some $\theta_G,\theta_B,\theta_C\ge0$ and  $a\in  [0,\eta_0)$.  Let
$u_0:\O\to \tX_{a}^0$ be strongly $\F_0$ measurable. Then the following
assertions hold.
\begin{enumerate}
\item There is a unique mild solution $U$ of \eqref{SEtype} with paths in
$C([0,T];\tX_a)$ a.s..

\item\label{en:mainthm2} For every $\delta,\lambda>0$ with
\[
\delta+a+\lambda <\min\{1-\theta_G,\tfrac12-\theta_B,
\tfrac12-\theta_C,\eta_0\}
\]
there exists a version of $U$ such that $U-P(\cdot,0)u_0$ in
$C^\lambda([0,T];\tX_{\delta+a})$ a.s..

\item\label{eq:mainthm3} If $\delta,\lambda>0$ are as in \ref{en:mainthm2} and
if $u_0\in \tX_{a+\delta+\lambda}$ a.s., then $U$ has a version with paths in
$C^\lambda([0,T];\tX_{\delta+a})$ a.s..
\end{enumerate}
\end{theorem}

\begin{proof}
Assertions (1) and (2) can be reduced to the case
\[
\left\{\begin{aligned}
dU(t) & = (A(t)U(t) + \tilde{F}(t,U(t)) + \tilde{B}(t,U(t))\,dW_{H}(t), \  t\in [0,T],\\
 U(0) & = u_0.
\end{aligned}
\right.
\]
taking $\tilde{F} = F +\Lambda_G G$ and $\tilde{B} = (B,\Lambda_C C)$ and
$H= H_1\times H_2$. The theorem now  follows from \cite[Theorem~6.3]{Vna}.
In view of (2), for assertion (3) we only have to show
 that $P(\cdot,0)u_0$ has the required regularity, which is proved in
 \cite[Lemma 2.3]{Vna}.
We note that, in order to apply the above results from \cite{Vna} here,
 one has to replace in \cite{Vna} the real interpolation spaces
of type $(\eta,2)$  by complex interpolation spaces of exponent $\eta$.
This can be done using the arguments given in \cite{Vna}.
\end{proof}

\section{Boundary noise\label{sec:ex1}}
 Let $S\subseteq \R^d$ be a bounded domain with $C^2$-boundary
 and outer unit normal vector of $n(s)$. On $S$ we consider the
stochastic equation with boundary noise
\begin{align}\label{eq:SE1}
\frac{\partial u}{\partial t}(t,s)
  &= \mathcal{A}(t,s,D) u(t,s)+ f(t,s,u(t,s)) \\
& \qquad  + b(t,s,u(t,s)) \, \frac{\partial w_1}{\partial t}(t,s),
 && s\in S, \ t\in (0,T], \notag\\
\mathcal{B}(t,s,D) u(t,s) &= {G}(t,u(t,\cdot))(s)+ \tilde{C}(t,u(t,\cdot))(s)
  \frac{\partial w_2}{\partial t}(t,s)\, , && s\in \partial S, \ t\in (0,T],
\notag \\
u(0,s) &= u_0(s), && s\in S. \notag
\end{align}

\noindent
Here $w_k$ are Brownian motions as specified below, and we use the
differential operators
\begin{align*}
\mathcal{A}(t,s,D) &= \sum_{i,j=1}^d D_i \big( a_{ij}(t,s) D_j \big) +
a_0(t,s),\quad
  \mathcal{B}(t,s,D) = \sum_{i,j=1}^d a_{ij}(t,s) n_i(s) D_j.
\end{align*}
For simplicity we only consider the case of a scalar equation, but systems
could be treated in the same way, cf.\ e.g.\ \cite{Amnonhom, DLS}.

\let\ALTERWERTA\theenumi
\let\ALTERWERTB\labelenumi
\def\theenumi{(A1)}
\def\labelenumi{(A1)}
\begin{enumerate}
\item \label{as:abcoef} We assume that the coefficients of $\mathcal{A}$
and $\mathcal{B}$ are real and satisfy
\begin{align*}
a_{ij} &\in C^{\mu}([0,T];C(\overline{S})), \ a_{ij}(t, \cdot)\in
C^1(\overline{S}), \ D_k a_{ij}\in C([0,T]\times\overline{S}),
\\ a_0&\in C^{\mu}([0,T],L^d(S))\cap C([0,T];C(\overline{S}))
\end{align*}
for a constant $\mu\in(\tfrac{1}{2},1]$ and all $i, j, k=1, \ldots, d$ and
$t\in [0,T]$. Further, let $(a_{ij})$ be symmetric and assume that there is a
$\kappa>0$ such that
\begin{equation}\label{unifell}
\sum_{i,j=1}^d a_{ij}(t,s) \xi_i \xi_j \ge \kappa |\xi|^2 \qquad
\text{for all \ } s \in
\overline{S},\: t \in [0,T],\: \xi \in \R^d.
\end{equation}

\end{enumerate}
\let\theenumi\ALTERWERTA
\let\labelenumi\ALTERWERTB

In the following we reformulate the problem \eqref{eq:SE1} as
\eqref{SEtype} thereby giving \eqref{eq:SE1} a precise sense.
Set $X = L^p(S)$ for some $p\in (1,\infty)$. Let $\alpha\in [0,2]$
satisfy $\alpha-\frac1p\neq 1$. We introduce the space
\[ H^{2,p}_{\mathcal{B}(t)}(S) =\begin{cases}
\big\{f\in H^{\alpha,p}(S):\ \mathcal{B}(t,\cdot,D) f = 0\big\},
     & \alpha-\frac1p>1,\\
 H^{\alpha,p}(S),& \alpha-\frac1p<1,
\end{cases}\]
where $H^{\alpha,p}(S)$ denotes the usual Bessel-potential space
(see \cite{Tr1}).  We also set
\[\tX_\eta = H^{2\eta,p}(S)\qquad \text{for all \ } \eta\geq 0.\]
We further define
$A(t):D(A(t))\to X$ by $A(t)x= \mathcal{A}(t,\cdot, D)x$ and
\[D(A(t)) = \{x\in H^{2,p}(S): \mathcal{B}(t,\cdot,D) x = 0\}
=H^{2,p}_{\mathcal{B}(t)}(S).\]

\begin{lemma}\label{lem:generationandHeta}
Let $X=L^p(S)$ and  $p\in (1, \infty)$. Assume that \ref{as:abcoef}
is satisfied. The following assertions hold.

\begin{enumerate}
\item\label{en:lemmapart1} The operators $A(t)$, $t\in [0,T],$ satisfy
(AT) and the graph norms of $A(t)$ are uniformly equivalent
with  $\|\cdot\|_{H^{2,p}(S)}$.  In particular,
$(A(t))_{t\in [0,T]}$ generates a unique strongly continuous evolution
family $(P(t,s))_{0\leq s\leq t\leq T}$ on $X$.

\item\label{en:lemmapart2} We have $X_\theta^t
=H^{2\theta,p}_{\mathcal{B}(t)}(S)$  for all $\theta\in (0,1)$ with
$2\theta - \frac1p \neq 1$, as well as  $X_\eta^t = \tX_\eta=H^{2\eta,p}(S)$
for all $\eta\in [0,\tfrac12+\frac{1}{2p})$, in the sense of isomorphic Banach
spaces. The norms of these isomorphisms are  bounded
uniformly for $t\in [0,T]$.

\item\label{en:lemmapart3} Let $p\in[2,\infty)$.
Then condition \ref{as:isom} holds with $\eta_0=1/2$.
\end{enumerate}
\end{lemma}

\begin{proof}
(1): \ See \cite{Acq} and \cite{Ya90}. Note that $A(t)^*$ on $L^{p'}(S)=X^*$
is given by $A^*(t)\varphi=\mathcal{A}(t,\cdot, D)\varphi$ with
$D(A(t)^*) = H^{2,p'}_{\mathcal{B}(t)}(S),$
 and thus also $(A(t)^*)_{0\le t\le T}$ satisfies (AT).

(2):\ Let $\theta\in (0,1)$ and $p\in (1, \infty)$ satisfy
$2\theta - \frac1p \neq 1$. Then Theorem~5.2 and Remark~5.3(c)
in  \cite{Amnonhom} show that
\begin{align}\label{interpolation}
X_{\theta}^t = [L^p(S), D(A(t))]_{\theta}
  = [L^p(S), H^{2,p}_{\mathcal{B}(t)}(S)]_{\theta}
=H^{2\theta,p}_{\mathcal{B}(t)}(S)
\end{align}
isomorphically, see also \cite[Theorem 4.1]{Se} and \cite[Theorem 1.15.3]{Tr1}.
Inspecting the proofs given in \cite{Se} one sees that the isomorphisms in
\eqref{interpolation} are bounded uniformly in $t\in [0,T]$.
Similarly, if $2\theta-\frac1p <1$, then
$X_{\theta}^t = H^{2\theta,p}_{\mathcal{B}(t)}(S) = H^{2\theta,p}(S)
= \tX_{\theta}.$

(3): This is clear from (1), (2) and the definitions. Note that the
spaces $\tX_{\eta}^t$ are UMD spaces with type $2$ because they are
isomorphic to closed subspaces of $L^p$-spaces with $p\in [2,
\infty)$.
\end{proof}

\begin{remark}
Let the constant $w\ge0$ be given by  (AT). In problem \eqref{eq:SE1} we replace
$\mathcal{A}$ and $f$ by $\mathcal{A}-w$ and $f+w$, respectively,
without changing the notation.
This modification does not affect the assumptions (A1) and (A2), and from now
 we can thus take $w=0$ in (AT).
\end{remark}

Next, we apply Theorem 9.2 and  Remark 9.3(e) of \cite{Amnonhom}
in order to construct the operators $\Lambda_C(t)$ and $\Lambda_D(t)$.
In \cite{Amnonhom} it is assumed that $\partial S\in C^\infty$.
However, the results from \cite{Amnonhom} used below remain valid
under our assumption that $\partial S\in C^2$, due  to  Remark~7.3 of
\cite{Amnonhom} combined with Theorem~2.3 of \cite{DDHPV}.

Let $t\in [0,T]$. In view of our main Theorem~\ref{thm:mainexdm}
we consider only $p\ge 2$ and $\alpha\in (1,1+\frac1p)$ though some of the
results stated below can be generalized to other exponents. Let
\[Y = \partial W^{\alpha,p}(S) := W^{\alpha-1-1/p,p}(\partial S)\]
be the Slobodeckii space of negative order on the boundary which is defined
via duality e.g.\ in (5.16) of  \cite{Amnonhom}. Let $y\in Y$.
Theorem 9.2 and  Remark 9.3(e) of \cite{Amnonhom} give
a unique weak solution $x\in H^{\alpha,p}(S)$ of the elliptic problem
\begin{align*}
\mathcal{A}(t,\cdot, D) x  &= 0 \qquad \text{on} \ S, \\
\mathcal{B}(t,\cdot,D) x &= y \qquad  \text{on} \ \partial S.
\end{align*}
(Weak solutions are defined by means of test functions
$v\in H^{2-\alpha,p'}(S)$, see \cite[(9.4)]{Amnonhom}.)
We set $N(t)y:=x$. Formula (9.15) of \cite{Amnonhom}  implies
that the `Neumann map'  $N(t)$ belongs
$\calL(\partial W^{\alpha,p}(S),H^{\alpha,p}(S))$ and that the
 map $N(\cdot):[0,T]\to \calL(\partial W^{\alpha,p}(S),H^{\alpha,p}(S))$
is continuous.

Concerning the other terms in the first line of \eqref{eq:SE1} and the
noise terms, we make the following hypotheses.

\let\ALTERWERTA\theenumi
\let\ALTERWERTB\labelenumi
\def\theenumi{(A2)}
\def\labelenumi{(A2)}
\begin{enumerate}
\item\label{as:fb} The functions $f,b:[0,T]\times\O\times S\times \R\to \R$ are
jointly measurable, adapted to $(\F_t)_{t\geq 0}$, and Lipschitz functions and
of linear growth in the third variable, uniformly in the other variables.
\end{enumerate}
\let\theenumi\ALTERWERTA
\let\labelenumi\ALTERWERTB

\let\ALTERWERTA\theenumi
\let\ALTERWERTB\labelenumi
\def\theenumi{(A3)}
\def\labelenumi{(A3)}
\begin{enumerate}
\item \label{as:noise} For $k=1,2$, the process $w_k$ can be written in the
form $i_k W_{H_k}$, where $i_1\in\g(H_1,L^r(S))$ for some $r\in
[1,\infty)$ and $i_2\in \g(H_2,L^s(\partial S))$ for some $s\in [1,
\infty)$, and $W_{H_1}$ and $W_{H_2}$ are independent $H_k$--cylindrical
Brownian motions with respect to $(\F_t)_{t\geq 0}$.
\end{enumerate}
\let\theenumi\ALTERWERTA
\let\labelenumi\ALTERWERTB

Supposing that \ref{as:fb} holds, we define $F:[0,T]\times\O\times X\to X$ by setting
$F(t,\omega,x)(s) = f(t,\omega, s, x(s))$. Then $F$ satisfies \ref{as:LipschitzFtype}.
We further define the function $B(t,\omega,x)h$ on $S$
for $(t,\omega,x)\in [0,T]\times\O \times X$ and $h\in H_1$ by means of
\begin{equation}\label{eq:Bh}
B(t,\omega,x) h = b(t,\omega,\cdot,x(\cdot)) \, i_1 h
\end{equation}
In Examples \ref{ex:bLr} and \ref{ex:white}
we give conditions on $w_1$ and $\theta_B$ such that $(-A)^{-\theta_B}B$
maps $[0,T]\times\O \times X$ into $\g(H_1,X)$ and
 \ref{as:LipschitzBtype} holds.

Assumption \ref{as:noise} has to be interpreted in the sense that
\[w_k(t,s) = \sum_{n\geq 1} (i_k h_n^k)(s) W_{H_k}(t) h_n^k,
\qquad  t\in \R_+, s\in S, \ k=1,2,\]
where $(h_n^k)_{n\geq 1}$ is an
orthonormal basis for $H_k$, and the sum converges in $L^r(S)$
 if $k=1$ and in $L^{s}(\partial S)$
if $k=2$. We note that then $(w_k(t,\cdot))_{t\ge0}$ is a Brownian motion
with values in $L^r(S)$ and $L^{s}(\partial S)$, respectively. Conversely,
if $(w_k(t,\cdot))_{t\ge0}$, $k=1, 2$, are independent Brownian motions with
values in $L^r(S)$ and $L^{s}(\partial S)$, then
we can always construct $i_k$ and $W_{H_k}$ as above, cf.\
Example~\ref{ex:i2q} below.

We recall that $H^{\alpha,p}(S)=X_\frac{\alpha}{2}^t$ for $t\in[0,T]$ and
$\alpha\in(1,1+\frac1p)$ by Lemma~\ref{lem:generationandHeta}(2). Moreover,
the operator $A(t)$ has bounded imaginary powers in $X$ (uniformly in $t\in
[0,T])$, see e.g.\ Example~4.7.3(d) and Section 4.7 in \cite{Am}. It then
follows that
 \begin{equation}\label{frac-compl}
 H^{\alpha,p}(S) =X_\frac{\alpha}{2}^t=D((-A(t))^\frac{\alpha}{2})
\end{equation}
with uniformly equivalent norms for $t\in [0,T]$, see e.g.\
\cite[Theorem~1.15.3]{Tr1}.
Therefore,  the extrapolated operator $A_{-1}(t)$ maps $H^{\alpha,p}(S)$
into $X_{\frac{\alpha}{2}-1}^t$, and hence
\[\Lambda(t) =\Lambda_G(t)= \Lambda_C(t):=- A_{-1}(t) N(t)\in
\calL(Y,X_{\frac{\alpha}{2}-1}^t)
\]
with uniformly bounded norms for $t\in [0,T]$. Let
$\theta\in[ 1-\frac{\alpha}{2},1]$. As above, we further obtain
$X_{\frac{\alpha}{2}-1+\theta}^t= H^{\alpha-2+2\theta,p}(S)
\hookrightarrow X$,  so that
\begin{equation}\label{eq:ALambdaC}
(-A(t))^{-\theta} \Lambda(t)\in
          \calL(Y, H^{\alpha-2+2\theta,p}(S))
\end{equation}
with uniformly bounded norms for $t\in [0,T]$.

In order to relate the boundary noise term in \eqref{eq:SE1} with
\eqref{SEtype}, we set
\[(C(t,\omega,x) h)(s) = \tilde{C}(t,\omega,x)(s) (i_2 h)(s)\]
for $h\in H_2$. We aim at the mapping property
$C(t,\omega,x):H_2\to Y=\partial W^{\alpha,p}(S)$ since it will enable us
to verify the hypothesis \ref{as:LipschitzCtype}. In fact,
if $C(t,\omega,x)\in\calL(H_2, Y)$
then $(-A(t))^{-\theta} \Lambda(t) C(t,\omega,x)$ maps $H_2$
continuously into $H^{\alpha-2+2\theta,p}(S)\hookrightarrow X$ if
$\theta\in[ 1-\frac{\alpha}{2},1]$. In Examples \ref{ex:i2p} and \ref{ex:i2q}
we give conditions on $\tilde{C}$ and $i_2$ implying \ref{as:LipschitzCtype}
for $C$.  The deterministic boundary term $G$ can be treated in a similar way.

We want to present a variational formulation of \eqref{eq:SE1}, starting
with an informal discussion. Let $\varphi\in \Gamma_t$,
where $\Gamma_t$ is given by \eqref{eq:Gammat}. Then $\varphi(r)\in
D(A(r)^*)=W^{2,p'}_{\mathcal{B}(r)}(S)$. Formally, multiplying
\eqref{eq:SE1} by $\varphi$, integrating over $[0,t]\times S$,  integrating
by parts and interchanging the order of integration, we obtain
 that, almost surely,
\begin{align}
\nonumber \int_S [u(t,s) \varphi(t)(s) - u_0(s)  \varphi(0)(s)] \, ds
   &=\int_0^t\! \int_{S} u(r,s) [\mathcal{A}(r,\cdot,D)\varphi(r)+\varphi'(r)](s)\,ds \,dr\\
   \label{eq:varmotivation}
   & \; + \!\int_0^t\!\int_S f(r,s,u(r,s)) \varphi(r)(s) \, d s \, dr \\ \nonumber
& \; +\! \int_S\! \int_0^t b(r,s,u(r,s)) \varphi(r)(s) \, d w_1(r,s) \, ds
+T_1 .
\end{align}
In the boundary term $T_1$ the part with $\nabla \varphi(r)$ disappears since
 $\varphi(r)\in D(A(r)^*)$, and the other term is given by
\begin{align*}
T_1 & = \int_{\partial S} \int_0^t \mathcal{B}(r,\cdot,D) u(r,\cdot)
\text{tr}(\varphi(r))\, dr \, d\sigma \\ & = \int_{\partial S}
\int_0^t {G}(r,u(r,\cdot)) \text{tr}(\varphi(r)) \, dr \, d\sigma  +
\int_{\partial S} \int_0^t \tilde{C}(r,u(r,\cdot))
\text{tr}(\varphi(r)) d w_2(r,\cdot) \, d\sigma
\end{align*}
where $\text{tr}$ denotes the trace operator on
$W^{2,p'}_{\mathcal{B}(r)}(S)$.

We now start from  the equation \eqref{eq:varmotivation} and rewrite it using
\eqref{weak-int} and the notation introduced above. Setting
$u(t,s) =: U(t)(s)$, equality \eqref{eq:varmotivation} becomes
\begin{align}\label{eq:varmotivation2}
\lb U(t), \varphi(t)\rb - \lb u_0, \varphi(0)\rb
&=  \int_0^t \lb U(r), (\mathcal{A}(r,\cdot,D)\varphi(r)+\varphi'(r) \rb \, dr
 + T_1\\
 &\;\; +\!\int_0^t\! \lb F(r,U(r)), \varphi(r) \rb \, dr
 + \int_0^t\! \! B(r,U(r))^* \varphi(r) \, d W_{H_1}(r), \notag
\end{align}
and the boundary term yields
\begin{align*}
T_1&= \int_0^t \lb {G}(r,U(r)), \text{tr}(\varphi(r))\rb  \, dr
+ \int_0^t C(r,U(r))^* \text{tr}(\varphi(r)) \, d W_{H_2}(r).
\end{align*}
Here the brackets denote the duality pairing on $L^p(S)$ and $L^p(\partial S)$,
respectively. We claim that for all $x\in W^{2,p'}_{\mathcal{B}(t)}(S)$ it holds
\[\text{tr}(x) = \Lambda(t)^*x = (-A_{-1}(t) N(t))^* x.\]
Indeed, let $\alpha\in (1, 1+\frac1p)$, $x\in W^{2,p'}_{\mathcal{B}(t)}(S)=D(A^*(t))$
 and  $y\in Y=\partial W^{\alpha,p}(S)$. Then we have  $N(t) y\in H^{\alpha,p}(S)$ and
 $ a(t)\nabla x \cdot n = 0$ on $\partial S$. Observe that $\Lambda(t)^*$ maps
 $D(A(t)^*)$ into $Y^*$.  Integrating by parts and using formula (9.4)
 of \cite{Amnonhom}, we obtain
\begin{align*}
\lb y, \Lambda(t)^*x\rb_Y & = \lb \Lambda(t) y, x \rb_{X_{-1}^t} = - \lb N(t)y, A(t)^* x \rb_X
\\ &  = - \int_S N(t) y [\nabla \cdot (a(t)\nabla x) + a_0(t) x] \, ds
\\ & = 0 + \int_S [(a(t) \nabla N(t) y) \cdot \nabla x + a_0(t)  (N(t) y) x] \, ds
= \lb y, \text{tr}(x)\rb_{Y},
\end{align*}
which proves the claim. Therefore, $T_1$ becomes
\begin{align*}
T_1 & = \int_0^t \lb \Lambda(t){G}(r,U(r)), \varphi(r)\rb  \, dr
                 +\int_0^t (\Lambda(t) C(r,U(r)))^* \varphi(r) \, d W_{H_2}(r).
\end{align*}
Combining this expression with \eqref{eq:varmotivation2} we arrive at the
definition of a variational solution to the stochastic evolution equation
\eqref{eq:SE1}, as introduced in Definition \ref{def:var}.
The above calculations thus motivate the following definitions. We say $u$ is a
{\em variational} (resp. {\em mild}) {\em solution} to \eqref{eq:SE1} if
$U(t)(s) = u(t,s)$ is a variational (resp. mild) solution to \eqref{SEtype}
with the above definitions of $A(t)$, $F$, $\Lambda_G$, $G$, $B$,  $\Lambda_C$,
$C$ and $W_{H_k}$. We can now state out main result.

\begin{theorem}\label{thm:mainexdm}
Let $p\in [2,\infty)$, $X=L^p(S)$, $\alpha\in (1,1+\tfrac1p)$,
$\theta_B\in [0,\frac12)$, $\theta_C\in (1-\frac{\alpha}{2},\frac12)$
and $\theta_G\in (1-\frac{\alpha}{2},1)$.
Assume that \ref{as:abcoef}--\ref{as:noise}  and
\ref{as:LipschitzGtype}--\ref{as:LipschitzCtype} hold,
where $C, G,B,\Lambda_C$ and $ \Lambda_G$ are defined above.
 Let $u_0:\O\to X$ be strongly $\F_0$-measurable.
Then the following assertions are true.
\begin{enumerate}
\item There exists a unique variational and mild solution $u$ of
\eqref{eq:SE1} with paths in $C([0,T];X)$ a.s..

\item For every $\delta,\lambda>0$ with $\delta+\lambda <\min\{1-\theta_G,
\frac12-\theta_B, \frac12-\theta_C\}$ there exists a version of $u$ such that
$u-P(\cdot,0)u_0$ in $C^\lambda([0,T];\tX_{\delta})$ a.s..

\item If $\delta,\lambda>0$ are as in  \ref{en:mainthm2} and
if $u_0\in \tX_{\delta+\lambda}$ a.s., then $u$ has a version with paths in
$C^\lambda([0,T];\tX_{\delta})$.
\end{enumerate}
\end{theorem}

Note that we need $\frac12-\theta_C<\frac{\alpha}{2} -
\frac12<\frac{1}{2p}$. Thus, if $\frac12 - \theta_B\geq \frac{1}{2p}$,
 $1- \theta_G\geq \frac{1}{2p}$ and the other assumptions in
Theorem~\ref{thm:mainexdm} hold, then we can take $\lambda,\delta\geq 0$
with $\delta+\lambda <\frac{1}{2p}$ and deduce that
 $u-P(\cdot,0)u_0$ belongs to $C^\lambda([0,T];\tX_{\delta})$ a.s..
If we also have $u_0\in H^{\frac1p,p}(S)$, then
we obtain a solution $u$ of \eqref{eq:SE1} with
paths in $C([0,T]; H^{2\delta,p}(S))$ for all $\delta<\frac{1}{2p}$.
In this case Sobolev's embedding (see \cite[Theorem 4.6.1]{Tr1}) implies that
\begin{align*}
u\in C([0,T]; L^q(S)) \ \ \text{for all} \ \ \left\{%
\begin{array}{ll}
    q<\frac{d p}{d-1}   & \hbox{if $d\geq 2$,} \\
    q<\infty, & \hbox{if $d=1$.} \\
\end{array}%
\right.
\end{align*}

\begin{proof}[Proof of Theorem \ref{thm:mainexdm}]
The existence and uniqueness of a mild solution with the asserted regularity
follows from Theorem~\ref{thm:mainexistencegeninitial} and the above observations.
The equivalence with
the variational solution is a consequence of Proposition \ref{prop:varmild}.
\end{proof}

We now discuss several examples under which \ref{as:LipschitzBtype} and
\ref{as:LipschitzCtype} hold. The hypothesis \ref{as:LipschitzGtype}
can be treated in the same way. We start with some observations concerning
Gaussian random variables $\xi$ with values in a Banach space $Z$,
see e.g.\ \cite{Bog}, \cite{BrzvNe} and the references therein. The
covariance $Q\in \calL(Z^*, Z)$
of $\xi$ is given by $Q x^* =\E (\lb \xi, x^*\rb \xi)$ for $x^*\in Z^*$. One
introduces an inner product $[\cdot,\cdot]$ on the range of $Q$ by setting
\begin{equation}\label{Q-sp}
[Q x^*,Qy^*] := \lb Q x^*, y^*\rb =\E(\lb \xi, x^*\rb \,\lb \xi, y^*\rb)
\end{equation}
for $x^*,y^*\in Z^*$,  and we define $\|Qx^*\|_H^2=[Qx^*,Qx^*]$.
The \emph{reproducing kernel Hilbert space} $H$ of $\xi$ is
 the completion of $Q Z^*$ with respect to  $\|\cdot\|_H$. Then
 the identity on $QZ^*$ can be extended to a continuous embedding
 $i:H\hookrightarrow E$,  and it holds $Q=ii^*$. On the other hand,
 the random variables $w_k(t,\cdot)$ in \ref{as:noise} are Gaussian
 with covariance $Q_k= t\,i_ki_k^*$ for all $t\ge 0$ and $k=1,2$.

\begin{example}\label{ex:i2p}
Let \ref{as:noise} hold with $H_2=L^2(\partial S)$.
Assume that covariance operator $Q_2\in \calL(L^2(\partial S))$ of $w_2$
is compact. Then there exist numbers $(\lambda_n)_{n\geq 1}$ in
$\R_+$ and an orthonormal system $(e_n)_{n\geq 1}$ in $L^2(\partial
S)$ such that
\begin{equation*}
Q_2= \sum_{n\geq 1} \lambda_n e_n \otimes e_n
\end{equation*}
Assume that
\begin{equation*}
\sum_{n\geq 1} \lambda_n \|e_n\|_{\infty}^2<\infty.
\end{equation*}
We observe that the operator $i_2$ is given by $i_2 = \sum_{n\geq 1}
\sqrt{\lambda_n} e_n \otimes e_n$ and belongs to $\calL(L^2(\partial S),
L^\infty(\partial S))$. Let $p\in[2,\infty)$. Assume that
$\tilde{C}:[0,T]\times\O\times L^p(S)\to L^{p}(\partial S)$ is
strongly measurable and adapted, as well as Lipschitz and of linear growth in
the third variable uniformly in $[0,T]\times\O$. Then
\ref{as:LipschitzCtype} holds for $C = \tilde{C}i_2$ with $a= 0$ and
every $\theta_C\in (1-\frac{\alpha}{2},\frac12)$, where $\alpha\in
(1,1+\frac1p)$.
\end{example}

\begin{proof}
Lemma~\ref{lem:sq-fc-Lp} implies that $i_2\in\gamma(H_2,L^p(\partial S))$.
Fix $t\in [0,T]$, $\omega\in \O$ and $x,y\in X=L^p(S)$. Denote $K =
\|i_2\|_{\calL(H_2, L^\infty(\partial S))}$. The embedding $L^p(\partial
S)\hookrightarrow Y = \partial W^{\alpha,p}(S)$ and \eqref{ideal} yield
\begin{equation}\label{eq:Csob}
\|C(t,\omega,x) - C(t,\omega,y) \|_{\g(H_2,Y)}
\lesssim_{p,\alpha} \|C(t,\omega,x) -
C(t,\omega,y)\|_{\g(H_2,L^p(\partial S))}.
\end{equation}
Furthermore, for $h\in H_2$ and $s\in S$ we have
\begin{align}\label{c-ctilde}
|((C(t,\omega,x) - C(t,\omega,y)) h)(s)| &=
|\tilde{C}(t,\omega,x)(s) - \tilde{C}(t,\omega,y)(s)| \, |i_2 h(s)| \notag\\
&\leq K |\tilde{C}(t,\omega,x)(s) - \tilde{C}(t,\omega,y)(s)|\, \|h\|_{H_2}.
\end{align}
Lemma \ref{lem:sq-fc-Lp} and the assumptions of the example
then imply that
\begin{align*}
\|C(t,\omega,x) - C(t,\omega,y)\|_{\g(H_2,L^p(\partial S))} &
\lesssim_p K \|\tilde{C}(t,\omega,x) -
\tilde{C}(t,\omega,y)\|_{L^p(\partial S)} \notag\\
&\leq  K L_{\tilde{C}} \|x-y\|_{L^p(S)}.
\end{align*}
Using \eqref{eq:ALambdaC}, we can now deduce
 the first part of \ref{as:LipschitzCtype}. The second
part is shown in a similar way.
\end{proof}

\begin{remark}
Note that in Example~\ref{ex:i2p} the noise could be a bit more irregular since in
\eqref{eq:Csob} one can still regain some integrability by choosing
$\alpha$ and $\theta_C$ appropriately.
\end{remark}

\begin{example}\label{ex:i2q}
Let $q\in (p, \infty]$ and $s\in [p,\infty)$ satisfy $\frac1p = \frac1q +
\frac1s$. We assume that $w_2$ is an $L^s(\partial S)$--valued Brownian
motion. Let $H_2$ be the reproducing kernel Hilbert space of the Gaussian
random variable $w_2(1)=w_2(1,\cdot)$ with covariance $Q$ and  $i_2$ be the
embedding of $H_2$ into $L^s(\partial S)$. Then  we have
$i_2 \in \g(H_2,L^s(\partial S))$  (cf.\ \cite{Bog}, \cite{BrzvNe} and
the references therein for details).  It is easy to check that
$t^{-1/2} w_2(t)$ also has the covariance $Q$ for $t>0$. Due to Proposition~2.6.1
in  \cite{KwWo} we thus obtain
\[ t^{-1/2} w_2(t) = \sum_{n\geq 1} \lb t^{-1/2} w_2(t), x_n^*\rb Q x_n^* \]
in $X$ a.s.\
for every orthonormal basis $(Q x^*_n)_{n \geq 1}$ of $H_2$. Therefore
\[ w_2(t) = \sum_{n\geq 1} \lb w_2(t), x_n^*\rb Q x_n^* \]
converges in $X$ a.s..
We now define $W_{H_2}(t): QL^{s'}(\partial S) \to L^2(\O)$  by setting
\[W_{H_2}(t) Q x^* = \sum_{n\geq 1} \lb w_2(t), x_n^*\rb \lb Q x_n^*, x^*\rb
  = \lb w_2(t),x^*\rb\]
for each $x^*\in L^{s'}(\partial S)$ and $t\ge0$.
Then we deduce $\|W_{H_2}(t) Q x^*\|_{2}^2 = \lb Q x^*, x^*\rb = \|Q x^*\|_{H_2}^2$
from \eqref{Q-sp}, and thus $W_{H_2}$ extends to a bounded operator from $H_2$
into $L^2(\O)$. It is easy to check that $W_{H_2}$ is the required cylindrical
Brownian motion with $w_2=i_2W_{H_2}$; i.e., \ref{as:noise} holds for $k=2$.
Assume
that $\tilde{C}:[0,T]\times\O\times X:\to L^{q}(\partial S)$ is strongly
measurable and adapted, as well as Lipschitz and of linear growth in the third
variable uniformly in $[0,T]\times\O$. Then \ref{as:LipschitzCtype} holds for
$C = \tilde{C}i_2$, where we take $a=0$, $\theta_C\in
(1-\frac{\alpha}{2},\frac12)$ and $\alpha\in(1,1+\frac{1}{p})$.
\end{example}
\begin{proof}
Fix $t\in [0,T]$, $\omega\in \O$ and $x,y\in X=L^p(S)$. We argue
as in the previous example, but in \eqref{c-ctilde} we consider
$\tilde{C}(t,\omega,x) - \tilde{C}(t,\omega,y)$ as an multiplication
operator from $L^s(\partial S)$ to $L^p(\partial S)$.
Using  H\"older's inequality and \eqref{ideal}, we thus obtain
\begin{align*}
\|C(t,\omega,x) - C(t,\omega,y)\|_{\g(H_2,Y)}
&\lesssim_{p,\alpha} \|C(t,\omega,x) -
C(t,\omega,y)\|_{\g(H_2,L^p(\partial S))}
\\ & \leq \|\tilde{C}(t,\omega,x) - \tilde{C}(t,\omega,y)\|_{L^{q}(\partial
S)} \|i_2\|_{\g(H_2,L^{s}(\partial S))}
\\ & \leq L_{\tilde{C}} \|x - y\|_{L^{p}(\partial
S)} \|i_2\|_{\g(H_2,L^{s}(\partial S))}.
\end{align*}
The first part of \ref{as:LipschitzCtype} now follows in view of
\eqref{eq:ALambdaC}. The second part can be proved in the same way.
\end{proof}

We now come to condition \ref{as:LipschitzBtype}.

\begin{example}\label{ex:bLr}
Assume that \ref{as:abcoef}--\ref{as:noise} hold with $r\in (d, \infty)$.
Then \ref{as:LipschitzBtype} is
satisfied for all $\theta_B\in (\frac{d}{2 r},\frac12)$.
\end{example}

\begin{proof}
Let $\frac1q= \frac1p+\frac1r$ and $\theta_B\in (\frac{d}{2 r},\frac12)$. As
in Example~5.5 of \cite{SchVplate} one can show
\begin{equation}\label{eq:sob1}
L^q(S)\hookrightarrow  X_{-\theta_B}^t,
\end{equation}
where the embedding is uniformly bounded for $t\in [0,T]$.
  Fix $t\in [0,T]$, $\omega\in\Omega$ and $x,y\in X=L^p(S)$.
Arguing as in the previous example, by means of
\eqref{eq:sob1}, \eqref{eq:Bh}, H\"older's inequality, (A2) and \eqref{ideal}
we can estimate
\begin{align*}
\|(-A(t))^{-\theta_B} &(B(t,\omega,x) - B(t,\omega,y))\|_{\g(H_1,X)}
\\ & \lesssim_{\theta_B, p,r,n} \|B(t,\omega,x) -B(t,\omega,y)\|_{\g(H_1,L^q(S)))}
\\ & \leq \|b(t,\omega,x) - b(t,\omega,y)\|_{L^p(S)} \|i_1\|_{\g(H_2, L^r(S))}
\\ & \leq L_b \|x - y\|_{L^p(S)} \|i_1\|_{\g(H_2, L^r(S))}.
\end{align*}
This proves the first part of \ref{as:LipschitzBtype}. The second
part is obtained in a similar way.
\end{proof}

Finally, we consider the white noise situation in the case $d=1$.
\begin{example}\label{ex:white}
Let $d=1$ and  $p> 2$ and assume that \ref{as:abcoef}--\ref{as:noise} hold
with $i_1=I$. Then \ref{as:LipschitzBtype} is satisfied
for all $\theta_B\in (\frac{1}{2p } + \frac{1}{4},\frac12)$.
\end{example}

\begin{proof}
Let $\frac1q = \frac1p + \frac12$ and
$\theta_B\in (\frac{1}{2p } + \frac{1}{4},\frac12)$.
Fix $t\in [0,T]$, $\omega\in \O$ and $x,y\in X=L^p(S)$.
Observe that $(-A(t))^{-\theta_B}$ can be extended to $L^q(S)$
where it coincides with the fractional power of the corresponding
realization $A_q(t)$ of $\mathcal{A}(t,\cdot,D)$ on $L^q(S)$
with the boundary condition $\B(t,\cdot,D)v=0$. We further obtain
\[D((-A_q(t))^{\theta_B})\hookrightarrow
  (L^q(S), H^{2,q}(S))_{\theta_B,\infty} \hookrightarrow
   [L^q(S), H^{2,q}(S)]_{\vartheta} = H^{2\vartheta,q}(S).\]
for $\vartheta\in (\frac{1}{2p } + \frac{1}{4},\theta_B)$ with
uniform embedding constants, see Sections 1.10.3 and 1.15.2 of \cite{Tr1}
and \eqref{interpolation}. Sobolev's embedding then yields that
$D((-A_q(t))^{\theta_B})\hookrightarrow
 C(\overline{S})$. Using also H\"older's inequality, we thus obtain
\begin{align*}
|[((-A(t))^{-\theta_B} \! (B(t,\omega,x) - B(t,\omega,y))h)](s)|
 & \lesssim_{\theta_B,p} \!\|(B(t,\omega,x) - B(t,\omega,y))h\|_{L^q(S)}
\\ & \leq \|b(t,\omega,x) - b(t,\omega,y)\|_{L^p(S)} \|h\|_{L^2(S)}
\\ & \leq L_b \|x - y\|_{L^p(S)} \|h\|_{L^2(S)}
\end{align*}
for all $s\in S$.
Now we can apply Lemma \ref{lem:sq-fc-Lp} to obtain that
\[\|(-A(t))^{-\theta_B} (B(t,\omega,x) -
B(t,\omega,y))\|_{\g(H_1,X)} \lesssim_{\theta_B,p,n}  L_b
\|x-y\|_{L^p(S)}.\]
The other condition  \ref{as:LipschitzBtype} can be verified in the same
way.
\end{proof}

In the next remark we explain why one cannot consider Dirichlet boundary
conditions
with the above methods. This problem was not stated clearly in \cite{Mas}.
In the one dimensional case with $S=\R_+$, a version of \eqref{eq:SE1}
with Dirichlet boundary conditions
has been treated in \cite{AB2} using completely other methods and
working on a weighted $L^p$ space on $\R_+$.
\begin{remark}\label{rem:neumann}
Since we are looking for a solution in $X=L^p(S)$, we have to require
that $\alpha-2+2\theta_C\geq 0$, see \eqref{eq:ALambdaC}. The restriction
$\theta_C<\frac12$ in Theorem~\ref{thm:mainexistencegeninitial} then leads to
$1-\frac{\alpha}{2}\leq \theta_C<\frac12$, so that $\alpha>1$. On the
other hand, in the case of Dirichlet boundary conditions one has $\partial
W^{\alpha,p}(S) = W^{\alpha-\frac1p,p}(\partial S)$ and the Neumann map
$N(t)$ has to be replaced by the Dirichlet map
$D(t)\in \B(\partial W^{\alpha,p}(S), W^{\alpha,p}(S))$, where
 $D(t) y := x\in W^{\alpha,p}(S)$ is the solution of the elliptic problem
\begin{align*}
\mathcal{A}(t,\cdot, D) x  &= 0 \qquad \text{on} \ S, \\
x &= y \qquad  \text{on} \ \partial S
\end{align*}
for a given $y\in \partial W^{\alpha,p}(S)$.
To achieve  that $\Lambda_C(t):=-A_{-1}(t) D(t)$ maps into $X_{-\theta_C}^t$,
we need that $H^{\alpha,p}(S) = H^{\alpha,p}_{\mathcal{B}(t)}(S)$, and
hence $\alpha-\frac1p<0$ in the Dirichlet case; which contradicts
$\alpha>1$ and $p\ge1$ .
\end{remark}

\def\cprime{$'$} \def\polhk#1{\setbox0=\hbox{#1}{\ooalign{\hidewidth
  \lower1.5ex\hbox{`}\hidewidth\crcr\unhbox0}}} \def\cprime{$'$}
\providecommand{\bysame}{\leavevmode\hbox to3em{\hrulefill}\thinspace}

\end{document}